\documentclass{article}
\usepackage{amsthm,amsfonts,amsbsy, amssymb,amsmath,graphicx}

\author{Vassily Olegovich Manturov}

\date{}

\title{Minimal diagrams of virtual links: II}

\newtheorem{thm}{Theorem}
\newtheorem{lm}{Lemma}

\newtheorem{re}{Remark}

\newcommand{\eps}{\varepsilon}

\begin{document}

\maketitle

\abstract{In the present paper we bring together minimality
conditions proposed in papers \cite{MaArx,MaArx2} and present some
new minimality conditions for classical and virtual knots and
links.}

\section{The main result}

This paper is a sequel of \cite{MaArx,MaArx2}. We deal with
virtual link diagrams and test whether these diagrams are minimal
with respect to the number of classical crossings. All necessary
definitions can be found in \cite{Ma,MaArx,MaArx2}.

Throughout the text, all virtual links are thought to be
non-split, unless otherwise specified. In any minimality theorem
for links we assume that the link in question has no split
diagrams. All virtual links are thought to be orientable in the
sense of atoms, see \cite{Ma,MaArx}.

First, let us formulate the two main theorems from \cite{MaArx}
and \cite{MaArx2}.

The main Theorem from \cite{MaArx} says the following:

\begin{thm} If a virtual link diagram $L'$ is $1$-complete and $2$-complete then
it is minimal. \end{thm}

The main Theorem from \cite{MaArx2} says the following:

\begin{thm}
Suppose the diagram  $K$ of a classical knot is good. Then it is
minimal in the classical category. In other words, if the diagram
$K$ has $n$ crossings then for any classical diagram representing
the same knot the number of crossing is at least $n$.
\end{thm}

In the sequel, we shall refer to these two theorems as ``The first
Theorem'' and ``The second Theorem''. The first theorem says that
if the span of the Kauffman polynomial for a virtual link diagram
is ``as large as it should be'' and the genus of the corresponding
atom is ``as large as it should be'' then the diagram is minimal.

The second theorem deals only with classical knots (not virtual
knots, and not links) and it has only one condition. This
condition says that no cell of the corresponding atom touches
itself at a crossing. This condition much stronger than the first
minimality condition used in the first theorem saying that {\em
the leading and the lowest term of the Kauffman state sum
expansion do not vanish}. The condition of the second Theorem says
that each of the two extreme coefficients equals precisely to one,
and we have only one non-trivial summand contributing to each of
them.

In fact, the condition of the second Theorem allows to consider
the cablings of the initial diagram $L$, the diagrams $D_{m}(L)$
where $m$ is a positive integer. It turns out that the main
condition of the second Theorem (that the diagram is good) is
hereditary: if it holds for $L$, then it should hold for
$D_{m}(L)$. This allows to establish minimality by passing to
cablings and some more tricks (see \cite{MaArx2}). Here we do not
worry about the thickness of the Khovanov homology \cite{Shu} of
the corresponding atom.

The first condition of the first Theorem (saying that the leading
and the lowest terms in the Kauffman state-sum expansion) are {\bf
not} hereditary: they may hold for $L$, but not for
$D_{m}(L),m\neq 1$, and usually they {\bf do not}. For instance,
they do not hold in the classical case. Thus, to establish the
minimality of a (virtual) knot or link diagram we have to handle
the genus thus adding one more condition on the Khovanov homology.

In the present paper, we wish to formulate a stronger statement
explaining the connection between the first Theorem and the second
Theorem, namely, we shall deal with the genus of links $D_{m}(L)$,
and $D_{m}(L{\#}{\bar L})$, where ${\bar L}$ is the mirror image
of the diagram $L$.

First, we deal with classical knots (to have the connected sum
operation well defined). At the end of the present paper, we shall
prove a theorem on virtual links.

If some auxiliary theorem or lemma admits a formulation for
virtual links, we formulate it for the general case though we
might need it only for the case of classical knots.

Suppose we have a virtual link diagram $L$ with $n$ classical
crossings such that the corresponding atom has genus $g$ and Euler
characteristic $\chi=2-2g$. Then the maximal possible span for the
Kauffman bracket of $L$ is estimated as (\cite{MaArx}):

\begin{equation}
span\langle L\rangle \le 4n+2 (\chi-2).
\end{equation}

Also, the Khovanov homology has thickness $T(L)$ defined by the
diagonals it lives between, see \cite{MaArx2}.

In \cite{MaArx2} (using a refined version of the result from
\cite{Weh}) we demonstrated that

\begin{equation}
T(L)\le 2+g(L).
\end{equation}

We shall refer to the corresponding equalities as ``the span of
the Kauffman bracket is as large as it should be'' and ``the
Khovanov homology is as thick as it should be''.

\newcommand{\ZZZ}{\zeta_{m}}

Let $K$ be a classical knot diagram with $n$ crossings, $N=2n$.
The usual estimate for $span\langle D_{m}\rangle$ is
$2(m^{2}+m)N+2m\chi(K{\#}{\bar K})-4$. It is so if the neither the
leading term nor the lowest term of the Kauffman state-sum
expansion vanishes. The following assymptotic theorem says that if
this length is smaller then expected, but assymptotically the
difference between the real length and the estimate is not very
large, then the initial diagram is minimal. Namely, we have

\begin{thm}[The assymptotic theorem]
Suppose that for some $\eps>0$ there is an infinite sequence
$i_{1}<i_{2}<\dots i_{m}<\dots$ of positive integers such that for
any positive integer $m$ we have

$span\langle D_{i_m}(K{\#}{\bar K})\rangle\ge
2({i_m}^{2}+{i_m})N+2{i_m}\chi(K{\#}{\bar
K})-4-(4-\eps)({i_{m}}^{2}+i_{m})$.

Then the diagram $K$ is minimal in the classical category.
\label{tmgy}
\end{thm}

We have the following

\begin{lm}
If a virtual link diagram $L$ is good then $T(L)\ge
g(L)$.\label{lmy}
\end{lm}

\begin{proof}
Indeed, one should just take the $A$-state with $v_{-}$ associated
to all circles and the $B$-state with $v_{+}$ associated to all
circles. The property that the diagram is good guarantees that the
first chain is a cycle, whence the second one is not a boundary.
Recalling the definition of the atom genus, we get the required
estimate.
\end{proof}

Now, having a good diagram $L$ of a classical knot, we see that
all diagrams $D_{m}(L{\#}{\bar L})$ are also good. By Lemma
\ref{lmy} we see that the diagram $L$ obviously satisfies the
condition of Theorem \ref{tmgy}. Here we may take a constant for
$\zeta_{m}$. Thus, the diagram $L$ is minimal.

So, the estimate for the thickness of the Khovanov homology for a
good diagram $L$ is that it is in between $2+g(L)$ and $g(L)$.

This immediately results in the following
\begin{thm}
Let $L$ be a good virtual link diagram with $n$ classical
crossings. Then any virtual link diagram $L'$ equivalent to $L$
has at least $n-2$ classical crossings.
\end{thm}

\begin{proof}
Indeed, suppose we have a diagram $L'$ equivalent to $L$ with $n'$
classical crossings. Then the thickness of $L'$ is at least
$g(L)$, so the genus $g(L')$ is at least $g(L)-2$, thus,
$\chi(L')\le \chi(L)+4$. From this we see that

\begin{equation}
4n+2(\chi(L)-2)=span \langle L\rangle= span\langle L'\rangle \le
4n'+2(\chi(L)+2),
\end{equation}

so $4(n-n')\le 8$, that completes the proof.
\end{proof}

We are still unable to use the trick with connected summation with
mirror image as in \cite{MaArx2} and prove the exact result (in
the unframed category). However, the reasonings with Khovanov
homology give a minimality estimate between $n-2$ and $n$
classical crossings.

\begin{proof}[Proof of Theorem \ref{tmgy}]
The proof goes in the same lines as that of the second Theorem
\cite{MaArx}.

Indeed, fix a positive integer $m$. Suppose there is a classical
diagram $K'$ having $n'$  crossings ($n'<n$) and representing the
same classical knot as $K$.

The diagrams $D_{i_m}(K{\#}{\bar K})$ and $D_{i_m}(K'{\#}{\bar
K}')$ generate isotopic knots. Denote $D_{i_m}(K{\#}{\bar K})$ by
$D_{m}$ and denote $D_{i_m}(K'{\#}{\bar K}')$ by $D'_{m}$. By
definition we have $\langle D_{m}\rangle=\langle D'_{m}\rangle$.
Also, set $\chi=\chi(K{\#}{\bar K}),\chi'=\chi(K'{\#}{\bar K'})$.

We have:

\begin{equation}
span \langle D_{m}\rangle \ge 4 {i_m}^{2}N+2(\chi_{m}-2) -
(4-\eps)(i_{m}^{2}+i_{m}),
\end{equation}
where $\chi_{m}=\chi(D_{i_m})$. The atom $V(D_{m})$ has
${i_m}^{2}N$ vertices,  $2{i_m}^{2}N$ edges and ${i_m}\Gamma$
$2$-cells, where $\Gamma=N+\chi$ is the number of the $2$-cells of
the atom $K{\#}{\bar K}$. Thus,

\begin{equation}
span \langle D_{i_m}\rangle \ge 2 ({i_m}^2+{i_m})N+2{i_m}
\chi-4-(4-\eps)(i_{m}^{2}+i_{m}).\label{ra1}
\end{equation}

For $D'_{m}$ we have:

\begin{equation}
span \langle D'_{m}\rangle \le 2 ({i_m}^2+{i_m})N+2{i_m}
\chi'-4.\label{ra1}
\end{equation}

Thus, taking into account $\langle D'_{m}\rangle=\langle
D_{m}\rangle$, we get

\begin{equation}2(i_{m}^{2}+ i_{m})(N-N')\le
(4-\eps)(i_{m}^{2}+i_{m})+2{i_m}(\chi'-\chi).
\end{equation}

This leads to a contradiction since $N-N'\ge 2$.
\end{proof}

\begin{re}
If we deal with classical link diagrams then for any link diagram
$L$ which is not good and for any positive integer $m>1$, the
leading coefficient in the Kauffman state-sum expansion for
$\langle D_{m}(L)\rangle$ is equal to zero. Thus, we can not apply
Theorem \label{tmyg} to the classical case directly.

In the virtual case, there are some diagrams $L$ which are not
good, but for every $m>1$, the state-sum expansion of the Kauffman
polynomial for $D_{m}(L)$ does not vanish. But here we can not
apply Theorem \ref{tmyg} directly because the connected summation
for virtual diagrams is not well defined.
\end{re}

Besides the asymptotic theorem (which works in the case of
classical knots only, because it uses the connected sum
construction), we also have the following generalisation of the
First theorem:

\begin{re}
Theorem \ref{tmgy}  works also in the case of {\em long} virtual
knots. The proof is literally the same because we have a
well-defined connect summation for such knots.

In the long virtual category, we have a lot of examples where the
first minimality condition is hereditary. Such examples obviously
give us minimal diagrams by theorem \ref{tmgy}.

Also, analogous theorems remain true for {\em tangles and braids}.
We shall discuss it in separate papers.
\end{re}


\begin{thebibliography}{100}

\bibitem{Kau} Kauffman, L.H. (1987), State Models and the Jones Polynomial,
{\em Topology}, {\bf 26} (1987), pp. 395--407.

\bibitem{KaV}
Kauffman, L. H. (1999), Virtual knot theory, {\em European Journal
of Combinatorics} {\bf 20}(7), pp. 662--690.

\bibitem{Jon} Jones, V. F. R. (1985), A polynomial invariant for links via Neumann algebras, {\em
Bull. Amer.  Math.  Soc.},  {\bf 129}, pp. 103--112.

\bibitem{Kh} Khovanov, M. (1997), A categorification of the Jones polynomial,
{\em Duke Math. J},{\bf 101} (3), pp.359-426.

\bibitem{Kup} Kuperberg, G. (2002), What is a Virtual Link?, www.arXiv.org,
math-GT$\slash$0208039, {\em Algebraic and Geometric Topology},
2003, {\bf 3}, 587-591.

\bibitem{Ma} Manturov, V.O., Teoriya uzlov, {\em Regular and Chaotic Dynamics},
Moscow-Izhevsk, 2005.

\bibitem{MaXi} Manturov V.O. (2003), Kauffman--like polynomial and
curves in $2$--surfaces, {\em Journal of Knot Theory and Its
Ramifications}, {\bf 12}, (8), pp.1145-1153.

\bibitem{MaLong} Manturov, V.O. (2004), Long virtual knots
and their invariants, {\em Journal of Knot Theory and Its
Ramifications}, {\bf 13} (8), pp. 1029-1039

\bibitem{MaKho} Manturov, V.O., The Khovanov Complex for Virtual
Links, Arxiv:GT/ 0501317.

\bibitem{MaArx} Manturov, V.O., Minimal diagrams of classical and virtual
links, ArXiv:GT/ 0501393.

\bibitem{MaArx2} Manturov, V.O., Minimal diagrams of classical knots, ArXiv:GT/ 0501510.

\bibitem{Shu} Shumakovitch, A., Torsion of the Khovanov homology, Arxiv:GT/
0405474.

\bibitem{Weh} Wehrli, S., A spanning tree model for the Khovanov
homology, Arxiv: GT/ 0409328



\end{thebibliography}
\end{document}